\documentclass[12pt]{article}
\usepackage{philhspstyle}

\def\cf{{\it cf.}}

\def\bv{{\bf v}}
\def\bn{{\bf n}}

\def\divv{\mathrm{div}}
\def\grad{\nabla}

\def\raw{\rightarrow}

\def\htop{h_{top}}
\def\laplacian{\Delta}
\def\br{{\bf r}}

\def\phit{{\phi_t}}

\def\reals{{\mathbb R}}
\def\R{{\mathbb R}}

\def\naturals{{\mathbb N}}

\def\ie{\textit{ie.}}
\def\cf{{\it cf.}}

\def\bv{{\bf v}}
\def\bn{{\bf n}}
\def\bu{{\bf u}}
\def\bx{{\bf x}}
\def\by{{\bf y}}

\def\raw{\rightarrow}

\def\br{{\bf r}}

\def\Inv{^{-1}}

\def\homeos{homeomorphisms}
\def\pA{pseudoAnosov}

\def\diffeos{diffeomorphisms}

\def\Lm{L^{met}}
\def\Lt{L^{top}}
\def\hOmega{\hat{\Omega}}
\def\homega{\hat{\omega}}

\begin{document}

\title{Exponential growth in two-dimensional \\topological fluid dynamics\footnote{For the proceedings of the IUTAM Symposium on
Topological Fluid Mechanics II, Cambridge, UK, July 2012}}

\author{Philip Boyland}

\maketitle

\begin{abstract}  
This paper describes   topological kinematics associated
with the stirring by rods of a two-dimensional fluid.
The main tool is the Thurston-Nielsen 
(TN) theory which implies that depending on the
stirring protocol the essential topological length of
material lines grows either exponentially or linearly.
We give an application 
to the growth of the gradient of a passively advected
scalar, the Helmholtz-Kelvin Theorem  then yields
applications to Euler flows. The main theorem shows that
 there are periodic 
stirring protocols for which generic initial vorticity yields a 
solution to Euler's equations which is not periodic and further, 
the $L^\infty$ and $L^1$-norms 
of the gradient of its vorticity grow exponentially in time. 
\end{abstract}

\section{Introduction}
Knots are an essential ingredient of three-dimensional
topological fluid dynamics. Their presence as flow or field
lines is a marker of 
 a certain level of complexity in the system. A fundamental
topological principle
 is that codimension-two is the knotting dimension,
so circles can be knotted in $3$D and two-spheres in $4$D. In
a two-dimensional flow,   points  are codimension two, but
can they be knotted?  Clearly not if we consider the static
problem, but if we consider their motion, it is plausible
that if they get sufficiently entangled during their  evolution
they can have implications for the surrounding flow field.

There is a substantial body of mathematics available to
understand this situation, most prominently, the Thurston-Nielsen
theory. This theory has many aspects and ramifications, but 
the part that is most directly applicable to two-dimensional
fluid dynamics concerns the growth of material lines. Specifically,
there are certain motions of points or rods in
a fluid which imply that the length of a  class of material 
lines must grow exponentially as they are advected by the fluid.
Further, the theory gives many methods to detect and/or construct
such motions as well as algorithms for computing the rate of
exponential growth.

In a two-dimesional incompressible fluid when a material
line is being stretched exponentially, then 
under the tangent map some  vectors are growing and 
others are shrinking exponentially. Heuristically,
this has implications for a (typical) passively advected scalar: 
its gradients
must also grow exponentially. For an Euler flow, the 
Helmholz-Kelvin Theorem says that the vorticity  
is passively advected and thus   its gradients
are growing exponentially. However, to make these arguments
precise, we need conditions on the advected scalar to make sure
it has nontrivial gradients and more importantly, we must eliminate
the possibility of 
 `` unfortunate coincidences'' where the gradients often
line up with a contracting direction in the fluid.

As is usual in Dynamical Systems Theory, 
these issues are dealt with by considering
only generic scalars. A set is considered topologically
large or generic if it is dense, $G_\delta$, 
\ie\ it is the intersection of a countable family
of open, dense sets. Thus with the appropriate
norm (and thus topology) on the set of  scalar functions
we only consider those in a carefully chosen dense,
 $G_\delta$-set. 
For the case of Euler flows, we use the  fact that
the initial vorticty determines the initial velocity field and
so the initial vorticities
 ``parameterize'' the collection of Euler flows. Our results
 concern the behaviour of typical values of this ``parameter''.

The second component of this paper  studies the situation
where a time-periodic stirring protocol gives rise to
a time-periodic fluid motion. This component does not
use the TN-theory but shares the two themes of growth
rate of material lines and generic hypothesis.
 For a time-periodic fluid motion the presence of a time-periodic
passively advected scalar implies that the flow is
integrable and thus, in particular, the growth rate
of advected material lines is at most linear. Once again,
the Helmholtz-Kelvin Theorem yields an immediate application
to Euler flows.

For more information on Thurston-Nielsen Theory 
see \cite{thurston} and \cite{FLP}, for its dynamical
systems applications  see \cite{bdamster}, and for its
fluid mechanical applications see \cite{BAS1} and \cite{thiffsurvey}.
For many more references (and pictures of fluids being
stirred by \pA\ protocols) see the papers by Thiffeault and Stremler
in this volume. %For background in Dynamical Systems
%theory see \cite{HK}. 
The theorems in this paper appear
in \cite{bdeuler} in a somewhat different form; see that
paper for details of the proofs.

\section{Topological kinematics}
\subsection{Basic definitions and terminology}
 We first formalize the stirring of a planar
body of fluid by moving rods.
The \textit{fluid region} is a smooth, one-parameter
family of smooth, multi-connected, compact, planar domains
denoted $M_t$. In this family 
the outer boundary is held fixed while the inner disks
move.  We shall always  assume time-periodicity
of the domains with period one, so $M_{t+1} = M_t$ for all $t$.
The moving inner regions are called the \textit{stirrers},
and they are perhaps permuted each period.

Since we will be using a fair amount of dynamical systems
terminology and  what is called a flow in dynamical systems is
called a steady flow in fluid mechanics, we shall adopt the
terminology  that 
a \textit{fluid motion} is a smooth one-parameter family of
diffeomorphisms, $\phi_t :M_0 \rightarrow M_t$, with $\phi_0 = id$.
We may
view $\phi_t$ as  Lagrangian fluid displacement map:
the particle at $\bx \in M_0$ at time $0$ is at 
$\phi_t(\bx) \in M_t$ at time $t$. However, at this point
we are making no assumption that the fluid 
 or its \textit{velocity field}  
\begin{equation*}
 \bu(\phi_t(\bx), t) := \frac{\partial \phi_t}{\partial t}(\bx)
\end{equation*} 
satisfies any particular equations.
Because the fluid does not penetrate
the moving boundary (expressed by $\phi_t :M_0 \rightarrow M_t$)
the velocity field
satisfies the boundary conditions $\bu\cdot\bn_i = 
\dot{B}_i  \cdot \bn_i$, with $B_i(t)$ the motion of
the $i^{th}$ boundary.  The fluid motion is \textit{incompressible
or area-preserving} if $\nabla\cdot\bu = 0$, or equivalently, 
the Jacobian of $\phi_t$ is identically one, $\det(D\phi_t)\equiv 1$.

We have assumed that the stirring protocol is
time-periodic. A special situation of importance
  is when the 
 velocity field is also time-periodic with the
same period. In this case we can study the time evolution using
the \textit{Poincar\'e or time-one map},
 $\phi_1:M_0\raw M_0 = M_1$, which
satisfies $\phi_n = \phi_1^n$, with the superscript indicating
repeated composition.

\subsection{One-dimensional growth rates}
A {material line} in the fluid is described
by  a smooth {arc} or {simple closed curve (scc)}
$\gamma$. The main focus here is how  
 the material line grows in length as it is
passively transported by the fluid motion, which means we are analyzing
 the length of $\phi_t\circ\gamma$
as a function of $t$. We will
use two different ways to measure the length of the curve,
one topological and the other metric.

\subsubsection{Metric growth}
Let $\ell_t(\gamma)$ denote the length of the curve
$\gamma$ with respect
to some smooth, periodic family of Riemannian metrics on the $M_t$. 
 The \textit{metric growth rate} of $\gamma$  
is the growth  of  
\begin{equation*}
\Lm_t(\gamma) := \frac{\ell_t(\phi_t\circ\gamma)}{\ell_0(\gamma)}
\end{equation*}
as a function of $t$.
There are many methods available to quantify growth
rate. The situation here is 
rather simple and we will usually just be bounding the growth
above and/or below by simple functions, namely, by $c \lambda^t$
with $\lambda > 1$ for \textit{exponential growth} and by
 $c\,t$ for \textit{linear growth}.

\subsubsection{Topological growth}\label{topogrow} 
The  \textit{topological one-dimensional growth rate}
is designed to depend on just the
topology of the stirring process and to be independent of various
details of the fluid motion. There are two main ideas.
First, we restrict consideration to so-called essential 
arcs and scc which  truly ``see the topology'' of the
evolving regions.  Second,  at each moment in time 
we don't compute the actual metric length of
the evolving curve  but rather we compute the shortest
length amongst curves with the same topology. Imagine
an advecting arc 
as  elastic; at each time we  will let it
 shrink back to the shortest length while maintaining
the endpoints on the same boundary circle. Maintaining
the curves topology means that in the shortening process
it is not allowed to pass through the solid stirrers.
 
These ideas are more formally defined  using homotopies.  
 Two arcs $\gamma_1,\gamma_2:[0,1]\raw X$
are \textit{homotopic} in $X$ if there exists a continuous
map $\Gamma:[0,1]\times [0,1]\raw X$ with $\Gamma(1,t) =
\gamma_1(t)$ and $\Gamma(2,t) = \gamma_2(t)$. Thus $\Gamma$
gives a continuous deformation from one arc to the other.

An \textit{essential arc} is required to  have its
endpoints on the boundary and two essential
arcs are homotopic if there is is homotopy 
$\Gamma$ between them with the additional property
that the homotopy keeps the endpoints on
the boundary, so $\Gamma(s, 0)$ and $\Gamma(s,1)$
are contained in the boundary for all $s$, Note that by continuity
of $\Gamma$   this means that $\gamma_1$ and $\gamma_2$
have their endpoints on the same pair of boundary circles.
In addition, for essential arcs $\gamma$ it is required that 
$\gamma$ is not homotopic to an arc wholly contained in
a single boundary circle. So, for example, an arc 
whose endpoints are on different boundary circles
is always essential.

We let $S^1$ denote the circle and so a scc in $X$ is
a continuous, injective map $\gamma:S^1\raw X$.   
Two scc $\gamma_1,\gamma_2: S^1\raw X$
are \textit{homotopic} in $X$ if there exists a continuous
map $\Gamma:[0,1]\times S^1\raw X$ with $\Gamma(1,t) =
\gamma_1(t)$ and $\Gamma(2,t) = \gamma_2(t)$.   
An \textit{essential scc} is one that is neither homotopic
to a point nor homotopic to a boundary circle.
Thus an essential scc   
encloses at least two but not all of the stirrers. All the
arcs and scc shown in Fig. 1.(b) are essential.

In both cases we denote the set of curves homotopic to
$\gamma$, \ie\ its homotopy class, as $[\gamma]$ 
and the topological length $\Lt(\gamma)$ of a curve will be the 
least length among curves in its homotopy class  
\begin{equation*}
\Lt(\gamma) := \min\{\ell(\sigma) \colon \sigma\in [\gamma] \}.
\end{equation*}
  The  \textit{topological
growth rate of the class} $[\gamma]$ is the growth of 
\begin{equation*}
\Lt_t(\gamma) := \frac{\Lt(\phi_t\circ\gamma)}{\Lt(\gamma)}. 
\end{equation*}
 Thus to compute the 
topological growth rate
we evolve the curve forward for time $t$ and then measure
the least length in its homotopy class. An crucial property  
 for   essential curves  $\gamma$ is that  
\begin{equation}\label{lowerbound}
\Lm_t(\gamma) \geq \Lt_t(\gamma).
\end{equation}

\subsection{Topological entropy}
A fundamental result  for two-dimensional iterated
$C^\infty$-\diffeos\ is that  the maximum exponential
metric growth rate of arcs is equal to the topological entropy
(\cite{newhouse1},  \cite{newhouse2}). Note the potential confusion with the current
terminology, namely, the  maximal \textit{metric}
one-dimensional growth is equal to the \textit{topological}
entropy. This entropy is so-called since it may be defined
using growth rates of distances using a topological metric, in contrast
to the measure-theoretic (or metric) entropy which requires
an invariant measure. The variational principle says that
the topological entropy is the supremum of the measure-theoretic
 entropies over all invariant measures: however, for
curves, the topological growth is a lower bound for the metric
growth. With all this potential confusion we hopefully clarify
matters by  emphasizing that our main concern here is the  
exponential growth of the lengths of one-dimensional curves.

 In analogy with the
result for iterated \diffeos, for a general 
two-dimensional fluid motion $\phi_t$ define its \textit{topological
entropy} as 
\begin{equation*}
\htop(\phi_t) = \sup\{\limsup_{t\raw\infty}
\frac{\log(\Lm_t(\gamma))}{t}\colon 
\gamma\ \text{is a smooth arc}\}.
\end{equation*}
In accord with equation \eqref{lowerbound}, if some essential
curve $\gamma$ has $\Lt_t(\gamma) \geq C\lambda^t$ for 
some constants $C>0$ and $\lambda > 1$, then
$\htop(\phi_t)\geq \log(\lambda) > 0$.

\subsection{Isotopy and braids}  
 As noted in Subsection~\ref{topogrow}, 
the topological growth rate of curves  depends only on the
rough topology of the stirrer motion. This motion, in turn,
determines what is called the isotopy class of the
map $\phi_1$.
 More precisely,  two \homeos\ $f, g:M_0\raw M_0$ are
\textit{isotopic} if there is a continuous family
of \homeos\ $h_t$ deforming one to the other with
$h_0 = f$ and $h_1 = g$. Given an essential
arc $\gamma$, applying the isotopy between $f$ and $g$ yields
a homotopy $H_t = h_t\circ \gamma$ between  
$ f\circ \gamma$ and  
$g\circ \gamma$. Thus we see that isotopic \homeos\
 give the same topological growth rate of an essential
curve $\gamma$.

The isotopy class of a stirring protocol can be visualized
and in fact characterized by the space-time trace of the
stirrers or  their \textit{braid}. The collection of braids 
forms a group and
the algebra of the braid gives a convenient
method to compute properties associated with the
trichotomy of the next subsection. 

Another  comment that will be relevant later
 is that  because our
stirring protocols are periodic, 
$\phi_1$ is isotopic to   $\phi_n$ for
all $n\in\naturals$, even when the velocity field
is not periodic.  

\subsection{The Thurston-Nielsen trichotomy}
 For the purposes of this paper, 
the \textit{Thurston-Nielsen theory} classifies stirred fluid motions
and their isotopy classes based upon the rate of their 
topological one-dimensional growth: is it exponential or
linear.  
\begin{theorem}[Thurston-Nielsen Trichotomy]
Let $M_t$ be periodic stirring protocol with
 fluid motion $\phi_t$.
  Then either 
\begin{enumerate}
\item {PseudoAnosov (pA):}  
there exist constants  $\lambda > 1$ {(the dilation)} and  
 $0 < C_1 < C_2$ such that  for every essential curve $\gamma$, 
\begin{equation*}
C_1\lambda^t \leq \Lt_t(\gamma)\leq C_2\lambda^t, 
\end{equation*} and thus $\htop(\phi_t)\geq \log(\lambda)>0$.
  \item {Finite order (fo):}
  there  exists a constant $K>0$ such that 
for every essential curve $\gamma$,
\begin{equation*} 
\Lt_t(\gamma)< K\, t.  
\end{equation*}
  \item {Reducible case:} (roughly stated)
  $M_0$ splits into $\phi_1$-invariant subsurfaces 
on which 1. or 2. holds.
 \end{enumerate} 
\end{theorem}
 We shall call a stirring protocol finite order,
\pA, or reducible according to its TN-type.
A critical feature of the TN-trichotomy is that in the finite
order and pA cases, {every} essential curve has the
same basic topological growth rate. The \pA\ case is the most
useful for applications because this rate is exponential.

It is important to emphasize that the TN theory only concerns
topological growth rates and thus it  
only gives  bounds for the more physical metric growth. Many
or perhaps most of the fluid motions arising from the 
same stirring protocol (\ie\ in the same isotopy class) could
have a greater growth rates. In particular, a finite
order stirring protocol of a real fluid will almost
certainly have regions of exponential metric growth while
its topological one-dimensional growth is always linear.
 
It is also important to note that the theorem was first proved
in the much broader context of isotopy classes of surface \homeos.
It is usually framed in terms of the existence of a special
map, the Thurston-Nielsen representative, which is present in
 each isotopy
class. When this special representative map is \pA, 
it has a wealth of nice
properties including a Markov partition which yields
a symbolic model of a mixing subshift of finite type. This
has many implications, for example, the \pA\ map 
is ergodic and mixing with respect to Lebesgue measure
and has a dense orbit and its set of periodic orbits is also
dense.  Handel's Isotopy Stability Theorem  shows that
most of the dynamical properties of the \pA\ map are 
present in any other map in the isotopy class, though they
could exist in a ``small'' invariant set (\cite{handel1}).
 
\section{Passive advection of scalars}
In this section we give two results concerning the 
behaviour of advected scalars under a stirring protocol.
The results continue to be strictly kinematic.
Given a fluid motion $\phi_t$, a function $\alpha:M_t\times\R 
\raw \R$ is called a  \textit{passively advected scalar} if 
it is constant on trajectories,  $\alpha_t(\phi_t(\bx)) =
\alpha_0(\bx)$, or equivalently, 
  $\partial \alpha_t(\phit(x))/\partial t = 0,$
where we have written $\alpha_t(\bx)$ for $\alpha(\bx, t)$.
In the language of global analysis one says that $\alpha_t$
is the \textit{push forward} of $\alpha_0$ and writes
$(\phi_t)_*(\alpha_0) = \alpha_t$, with $(\phi_t)_*(\alpha_0) =
\alpha_0\circ (\phi_t)\Inv$. 

Now for any function
$f:M_0\raw\R$ we may obtain a passively advected scalar simply
by defining $\alpha_t :=  (\phi_t)_*(f)$, and so in this sense
only the initial configuration  $\alpha_0$ matters.
For example, 
if $\alpha_0$ represents the initial concentration of dye
in a fluid, $\alpha_t =  (\phi_t)_*(\alpha_0)$ is the
concentration after time $t$.
On the other hand, in a physical fluid there may be scalar quantities
of importance which at each time $t$
are computed from the velocity field. 
In this case the advected scalar represents a conserved
quantity having physical meaning.  For example, in two
dimensions the curl,  $\omega_t = \nabla\times\bu$, is
a scalar which is passively advected in an Euler flow.

\subsection{Consequences of a \pA\ protocol}\label{pAcons}
For a fluid motion determined by a  \pA\ protocol
we know from the TN-trichotomy and equation
 \eqref{lowerbound} that the metric
length of essential curves is growing exponentially
fast. This means that tangent vectors to these curves
must be growing in length exponentially under the
action of the space derivative of the fluid motion 
$D\phi_t$. This implies that an eigenvalue of $D\phi_t$
is growing exponentially.
 For a   passively advected scalar
$\alpha_t$, since $\alpha_t = \alpha_0\circ (\phi_t)\Inv$, 
we have 
$\nabla \alpha_t = \nabla\alpha_0 (D\phi_t)\Inv.$
If the fluid motion is incompressible then $\det(D\phi_t) = 1$
and so  $(D\phi_t)\Inv$ also has an eigenvalue growing exponentially.
Thus as long as there is not a unfortunate coincidence where
$\nabla\alpha_0$ stays aligned with the stable eigen-direction
of $(D\phi_t)\Inv$, we have that $| \nabla \alpha_t |$ is
growing exponentially.

\begin{figure}
\centerline{
\includegraphics[width= .65\linewidth]{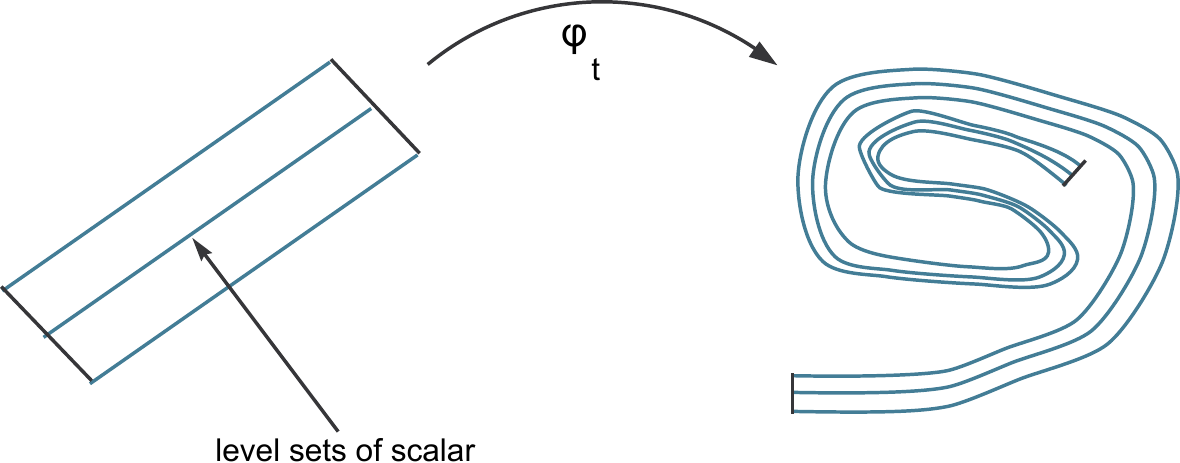}\hspace*{8mm}
\includegraphics[width= .30\linewidth]{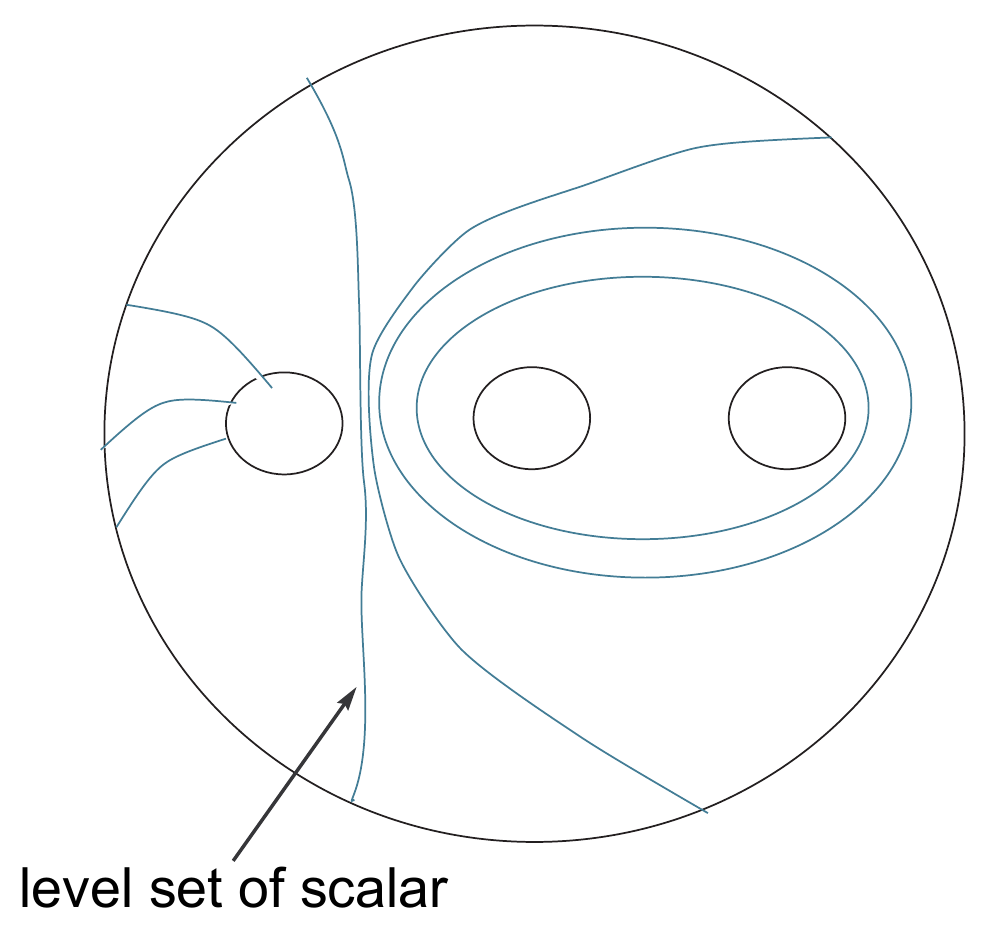}}
\vskip6pt
\caption{(a) Proof of Theorem 3.1 \qquad \qquad \qquad
(b) Proof of Theorem 3.2}
\end{figure}

There are a number of issues involved with making this
arguement rigorous;  the locations in 
  $\gamma$ where there is tangential growth will
be shifting in time and the unfortunate coincidence
could actually happen, especially if one is working
in a restricted class of physical interest.
One way around these difficulties
is to take a more geometric global viewpoint.
\begin{theorem}\label{expgrowthscalar}
  Let $M_t$ be a time-periodic stirring
protocol of {pA} type with incompressible fluid motion $\phi_t$.
If $\alpha_t$ is a passively advected scalar such that
its intial state   
$\alpha_0$ is a generic $C^2$-function, then
there are  positive constants $c, c'$  so that
\begin{equation*}
\sup_{\bx\in M_0} |\grad \alpha_t (\bx) | \geq c \lambda^t
\ \ \ \text{and}\ \ \ 
\int_{M_t} |\grad \alpha_t (\bx) | \geq c' \lambda^t
\end{equation*} 
for all $t \in \reals$, where   $\lambda>1$  is the dilation
of the \pA\ protocol.
\end{theorem} 
Here are  the main ideas in the proof.
First,  find a $C^2$-open, dense set $\mathcal{G}$ 
inside the Morse functions (functions with nondegenerate critical
points)
on $M_0$ so that $\alpha_0\in \mathcal{G}$ implies that $\alpha_0$
has a band of regular inverse images which are essential arcs or circles.
 The pA protocol forces a
stretch in length by $\lambda^t$. This coupled with area preservation
and transport of vorticity 
force the level sets of $\alpha_t = (\varphi_t)_\ast(\alpha_0)$
 to bunch up
in the transverse direction,
which  causes $\|\grad \alpha_t\|_\infty \rightarrow \infty$ like 
$\lambda^t$. See Fig. 1.(a)

\subsection{Time-periodic fluid motions}\label{scalarlinsec}
Now consider the special situation when the
velocity field of the fluid motion $\phi_t$ is time-periodic
with period  one 
and so $\phi_1$ is a Poincar\'e map.
In addition, we assume the existence of a passively advected scalar
that is also time-periodic. This again is rather special,
but if the advected scalar depends on the velocity field
(\textit{eg} its curl), it is automatically time-periodic when the
velocity field is. In this situation we show that if the scalar
field is typical in  the appropriate sense, then there is at 
most linear growth of the metric length of material lines.

 The main
observation required is that given a time-periodic
 passively transported scalar $\alpha_t = \alpha_{t+1}$ and
the Poincar\'e map $\phi_1$, the initial configuration  
  $\alpha_0$ satisfies for every $\bx$
\begin{equation*}
\alpha_0(\bx) = \alpha_1(\phi_1(\bx)) = \alpha_0(\phi_1(\bx)).
\end{equation*}
Thus  $\alpha_0$ is what is called
 an \textit{integral of motion} for the map $\phi_1$ and so
as long as $\alpha_0$ is sufficiently nondegenerate, its existence
precludes chaotic dynamics.
\begin{theorem}\label{lingrowthscalar}
Let $M_t$ be a time-periodic stirring
protocol  with  time-periodic  fluid motion $\phi_t$
and time-periodic transported scalar $\alpha_t$.
If the  initial  state of the scalar  
$\alpha_0$ has finitely many critical points
(for example, is $C^\omega$ and nonconstant
 or is $C^2$-generic), there exists
a constant $K>0$ so that 
\begin{equation*}
\Lm_t(\gamma) \leq K t
\end{equation*}
for all arcs and scc $\gamma$. Thus 
the {metric} one-dimensional growth is linear, and the Poincar\'e map
$\phi_1$ has zero topological entropy and zero 
Lyapunov exponents almost everywhere.
\end{theorem}

The geometric heart of the proof is this: the 
level sets of $\alpha_0$ must be preserved
by $\phi_1$ and for a generic function $\alpha_0$ the contours
$\alpha_0(\bx) =c$  are 
smooths arcs and circle for typical $c\in\R$. 
Thus  the dynamics of $\phi_1$
 consists of one-dimensional invariant subsets. Injective
maps of arcs and circles yield very simple dynamics
with    no dynamical entropy and using this 
in the entire fluid domain we have only
linear metric growth of arcs.

\subsection{Combining the two results}
Putting these two kinematic results together 
yields an  observation of relevance to Euler flows.
If a time-periodic fluid motion has 
a generic advected scalar field that depends
on the velocity field and thus is also time-periodic
we know from Theorem~\ref{lingrowthscalar} that there is 
linear metric one-dim growth of all curves (material lines).
On the other hand, we know from Theorem~\ref{expgrowthscalar} 
that a \pA\
stirring protocol always causes exponential stretching
of some material lines. These clearly can't coexist and so
if a fluid motion is stirred by a pA protocol
and  has a passively advected  generic scalar then
either the scalar is not dependent directly on the velocity field
or if it is, then the fluid motion is not time-periodic. 

\section{Euler fluid motions}

\subsection{Basic definitions and results} 
Now  we assume that the velocity field $\bu(\bx,t)$ of
the fluid motion $\phi_t$ satisfies the incompressible, constant
density ($\rho \equiv 1$), 
Euler equation
\begin{equation*}
\frac{D\bu}{Dt} =  - \grad p_t,   \ \ \ \ \
\divv(\bu) = 0,
\end{equation*}
with slip boundary conditions on the moving boundary.
In this case $\phi_t$ is called an \textit{Euler fluid motion}.
Kozonoi has shown that for the class of two-dimensional
problems with moving boundary  considered here
there are classical solutions (\cite{koz}). 
The results in this paper are all predicated on
the assumption that there is   
  a global strong solution with the
regularity of the initial data: our goal is to analyze its dynamics.

Recall that for two-dimensional,
divergence-free velocity fields
a classical result says 
that the {curl} coupled with the  {circulations}
around the boundary components and the {boundary
conditions} $\bu\cdot\bn_i = 
\dot{B}_i  \cdot \bn_i$ determine the field
completely. In a certain sense we will  view the collection
of Euler fluid motions as being ``parameterized'' by possible
initial vorticities and all our results concern dynamics
associated with typical values of these ``parameters''.  
 
The Helmholtz-Kelvin Theorem (1890's) is perhaps the
most important result about two-dimensional 
Euler fluid motions:  an incompressible
 fluid motion is Euler
if and only if its vorticity is passively transported
and circulations around all smooth
simple closed curves $C$ are preserved, or
\begin{equation*}
\frac{d}{dt}\oint_{\phi_t(C_i)} \bu \cdot d\br = 0
\end{equation*}
for each boundary circle $C_i$.
  The preservation of boundary circulation is a
necessary feature in multi-connected domains.
 
\subsection{An Exponential Growth Theorem}
Now we  consider an Euler fluid motion under the influence of a 
\pA\ stirring protocol. Theorem~\ref{expgrowthscalar}
 used with the Helmholtz-Kelvin
Theorem immediately yields:     
\begin{theorem}\label{expgrowththm}
  Let $M_t$ be a time-periodic stirring
protocol of {pA} type with Euler fluid motion $\phi_t$.
If the  initial vorticity 
$\omega_0$ is a generic $C^2$-function, then
there are  positive constants $c, c'$  so that
\begin{equation}\label{yudgrowth}
\sup_{\bx\in M_0} |\grad \omega_t (\bx) | \geq c \lambda^t
\ \ \ \text{and}\ \ \ 
\int_{M_t} |\grad \omega_t (\bx) | \geq c' \lambda^t
\end{equation} 
for all $t \in \reals$, where   $\lambda>1$  is the dilation
of the pA protocol.
 Thus $\|\laplacian \bu(\bx, t)\|_\infty =\|\grad \omega_t \|_\infty
\rightarrow \infty$ and $\|\bu_t\|_{C^2}\rightarrow\infty$, all like 
$\lambda^t$.
\end{theorem} 
There are numerous results in the stability literature 
concerning growth of $|\grad \omega_t|$ 
for perturbations of 
two-dimensional steady Euler fluid motions.    
 Yudovich (\cite{yud1}, \cite{yud2}) and others have shown linear
growth   and 
 Arnol'd (\cite{arnoldhyper}), Friedlander and  Vishik 
((\cite{friedlandervishik}) 
and others have shown the importance of exponential
growth.   
 The same basic mechanism is in play here
 as in those results, namely, the growth of
$\nabla \omega$ under the influence of the 
eigenvectors of the tangent map $D\phi_t$ as discussed
at the beginning of Subsection~\ref{pAcons}.

\subsection{Speculations on applications to more general Euler flows}
We   make a few very speculative remarks on how the exponential
growth induced by a \pA\  stirring protocols might
have implications for general Euler flow. In particular, it
 provides a topological
perspective on a version of the 
\textit{Yudovich Hypothesis/Conjecture}:
for generic initial vorticity a two-dimensional Euler fluid motion
satisfies the exponential growth in equation~\eqref{yudgrowth}
(\cf\ \cite{yudschnirl}).

There are two main ingredients in this topological perspective.
The first is that the exponential growth
of material lines caused by the stirrers could
just as well be achieved by  advecting points in 
the fluid (\cite{bowen}, \cite{ghost}). 
The TN-trichotomy as stated above requires   these  
 ``virtual stirrers'' or ``ghost rods'' to be periodic
points. However, it is clear that an aperiodic
``tangle'' of points evolving can also force exponential
 topological growth of a class of material lines.
The second ingredient would say
 that for typical initial vorticity a two-dimensional Euler
fluid motion always has such a collection of fluid
trajectories. If both these ingredients 
were known, virtually the same proof as Theorem~\ref{expgrowththm}
would yield the version of the Yudovich Conjecture just given.
While both ingredients seem quite reasonable, rigorous results
seem a long way off. The topological ingredient is probably
tractable but determining even a part of the dynamical evolution
of a typical two-dimensional Euler fluid motion  requires tools 
still to be developed.

 \subsection{A linear metric growth theorem}  
Because 
the vorticity $\omega_t$ is passively transported 
scalar  for an Euler fluid motion  as in Subsection~\ref{scalarlinsec}
 it follows that
when the fluid motion is time-periodic, the initial
configuration $\omega_0$ is an 
integral of motion for the Poincar\'e map $\phi_1$.
 Brown and Samelson (\cite{brownsam}) observed  using a theorem of Moser
 that if  the integral of motion $\omega_0$ is real analytic, then
 $\phi_1$ can't have a transverse homoclinic intersection.
 More generally using Theorem~\ref{lingrowthscalar}
 and the Helmholtz-Kelvin
Theorem we have:   
\begin{theorem}\label{lineargrowthcor}
Let $M_t$ be a time-periodic stirring
protocol  with  time-periodic  Euler fluid motion $\phi_t$.
If the  initial vorticity 
$\omega_0$ has finitely many critical points
(for example, is $C^\omega$ and nonconstant
 or is $C^2$-generic), there exists
a constant $K>0$ so that 
\begin{equation*}
\Lm_t(\gamma) \leq K t
\end{equation*}
for all arcs and scc $\gamma$. Thus 
the {metric} one-dimensional growth is linear, and the Poincar\'e map
$\phi_1$ has zero topological entropy and zero 
Lyapunov exponents almost everywhere.
\end{theorem}

The crucial
features of the theorem are the time-periodicity of the
Euler fluid motion and its two-dimensionality. Thus a
similar result   holds in any two-dimensional fluid region. 
It is also important to note what 
the theorem does not say. It does not say that the generic
initial vorticity gives rise to a time-periodic Euler flow,
and it does not say that among initial vorticities
which do give rise to time-periodic Euler fluid motions
linear growth is typical. Understanding these issues requires
a knowledge of the size of the set of 
  time-periodic Euler fluid motions which are not steady.
This is ill-understood at present. There
are reasons to believe that, in fact,  generic initial vorticity
never gives rise to a time-periodic Euler fluid motion and
Theorem~\ref{lineargrowthcor} could be useful in proving this.

\subsection{A dichotomy for time-periodic 2D Euler flows}
Theorem~\ref{lineargrowthcor} is similar to and inspired by
Arnol'd's result (\cite{arnoldeuler}) on 
steady, three dimensional $C^\omega$-Euler  flows which
we frame as a dichotomy:
\begin{enumerate}
\item For \textit{typical} vorticity the preservation
of the Bernoulli function $p + \|\bu\|^2/2$ implies that the flow is
integrable and typical orbits are confined to invariant tori
and annuli  and  the flow 
thus has zero entropy and zero Lyapunov exponents
almost everywhere.
 \item For the \textit{atypical} case of a Betrami
flow where the curl is parallel to
the velocity field, one may have measures with positive entropy,
positive Lyapunov exponents, etc.
\end{enumerate}

For \textit{time-periodic, two-dimensional} Euler fluid
we know from Theorem~\ref{lineargrowthcor} that in
analogy with 1., 
for {typical} initial vorticity there is integrability and
linear metric one-dimensional growth. For an analog
of the atypical case 2., perhaps the simplest
situation is when
the vorticity  is constant everywhere.
In this case one always has {time-periodic}
solutions for {every} time-periodic 
stirring protocols using classical potential
theory: if 
$C$ is a given constant and 
$(\Gamma_0, \Gamma_1,  \; \dots\;, \Gamma_m)
\in \reals^{m+1}$ with $\sum \Gamma_i = 0 $ is
a vector of circulations, then 
there exists a unique time-periodic Euler
 fluid motion $(M_t, \bu_t)$ with
$\omega_t  \equiv C$ for all time $t$ and
$\oint_{B_{i}(t)} \bu \cdot d\br = \Gamma_i$, 
for $i = 0, \dots, m$. The proof is essentially to  
solve the Poisson  equation $\triangle \Psi = C$
for the needed stream function $\Psi$ at each time $t$
and these stream functions depend on time-periodic data
and so are time-periodic as are the velocity
fields $\grad^\perp \Psi$. 

Now if we assume the stirring protocol is \pA\
with dilation $\lambda$, the resulting
fluid motion satisfies 
$\Lm(\gamma) \geq \Lt(\gamma)\geq c\lambda^t$
for any essential curve $\gamma$.
This implies that  the Poincar\'e map
$\phi_1$ has $\htop(\phi_1) \geq
\log(\lambda) > 0$, 
ergodic invariant measures with positive metric entropy 
and thus positive Lyapunov exponents
and more $\dots$.
Thus one dichotomy for time-periodic solutions to
Euler's equation in a time-periodic moving domain is:
\begin{enumerate}
\item For \textit{typical} initial vorticity 
there is integrability and
the linear growth theorem above.
\item For the \textit{atypical} case
of constant vorticity classical potential theory yields
time-periodic Euler fluid motions and with \pA\ protocols
these  have chaotic dynamics.
\end{enumerate}
However, once again we  emphasize that the
dichotomy is strictly conditional. It says nothing
about the existence or prevalence of non-steady,
time-periodic Euler fluid motions.

\subsection{An energy bound}
Conservation of energy is a fundamental property of
Euler fluid motions in stationary domains.
For   stirring protocols $M_t$   a simple
computation yields that the total energy 
$E = \frac{1}{2}\|\bu\|_2^2$
 is evolving so that
\begin{equation*}
\frac{dE}{dt} = 
- \sum \oint_{\phi_t(C_i)} p\, \dot{B_i}\cdot d\bn_i.
\end{equation*}
This reflects the fact that the fluid can do work
on the stirrers and vice versa and so the long term
behaviour of the energy in a stirred Euler fluid is
initially unclear especially in light of the possible exponential
growth of $\|\grad\omega_t\|_1$. 
We sketch an argument whose main ideas are due to
Steve Childress which shows that for time-periodic
stirring protocols the energy is
uniformly bounded as $t\raw\infty$.

Assume now that the Euler fluid motion $\bu$ has smooth initial
vorticity $\omega_0$ and circulation around the $i^{th}$ boundary
equal to $\Gamma_i$. By Helmholtz-Kelvin for all
 times $t$, $\bu(\bx, t)$ 
has the same circulations and further, $|\omega_t| \leq
\Omega := \max(|\omega_0|)$. 
 We now assume for simplicity that  the outer boundary
of $M_t$ is the unit circle $S^1$ and so $M_t$ is contained in the
unit disk $D$. We     further assume that each stirrer is circular
with area $a$. Now fix a time $t$ and suppress dependence on it.
We decompose $\bu$ into pieces whose energy will be bounded separately.

First, extend $\omega$ to a function $\homega$ on $\R^2$ with 
$\homega\equiv 0$ outside $D$ and $\homega\equiv \Gamma_i/a$
inside the $i^{th}$ stirrer. Using Biot-Savart we define
\begin{equation}
\bv(\bx) = \frac{1}{2\pi} \iint_{D}
 \frac{\homega(\by) (\by-\bx)^\perp}{|\by-\bx|^2}\; d\by.
\end{equation}
Thus $ \bv$ has vorticity = $\omega$ in $M$,   has circulations
$\Gamma_i$, and  $|\bv|\leq   \hOmega$ in $M$ where
$\hOmega = \max\{\Omega, |\Gamma_i|\}$  using 
\begin{equation}  \label{intbound} 
  \max_{\bx\in M}\{  \frac{1}{2\pi} \iint_{D}  
 \frac{1}  {|\by-\bx|}\; d\by\} = 1.
\end{equation} 
This yields $\|\bv\|^2 \leq \pi   \hOmega^2$.
    
Next, let $h$ be the harmonic function on $M$ with
$\partial h/\partial\bn = -\bv\cdot \bn_i$ on each boundary
circle $C_i$ and $\grad h$ 
has zero circulation around all the boundaries.
Finally, let $g$ be the harmonic function on $M$ with
$\partial g/\partial\bn = \dot{B_i}$ on each boundary
circle $C_i$ and $\grad g$ has zero circulation around
 all the boundaries.

We now have that $\bv + \grad h + \grad g$ has vorticity
$\omega$, circulations equal to $\Gamma_i$, and 
the normal component of its velocity on the
boundary circle $C_i$ is $\dot{B_i}\cdot\bn_i$ and so
it is   equal to $\bu$.
Because $(\bv + \grad h)\cdot \bn_i = 0$
on every boundary circle $C_i$ and is divergence free we have that 
$(\bv + \grad h)\perp \grad h $ and 
$(\bv + \grad h)\perp \grad g $ working in $L^2$. Thus twice
the energy of $\bu$ is
\begin{equation*}
\|\bu\|^2 = \|\bv + \grad h + \grad g\|^2
= \|\bv + \grad h\|^2 + \|\grad g\|^2
= \|\bv\|^2 - \|\grad h\|^2 + \|\grad g\|^2
\leq \|\bv\|^2  + \|\grad g\|^2
\end{equation*}
Now we reintroduce the time dependence and do the 
decomposition of $\bu(\bx, t)$ into three parts
at each time $t$. By Helmholz-Kelvin the circulations for
$\bu(\bx, t)$ and thus $\bv(\bx, t)$ are the same for all $t$
as is the value of $\hOmega$. Further,  the bound using
equation \eqref{intbound} is independent of $t$ and so 
  for all $t$, $\|\bv(\bx, t)\|^2 \leq \pi   \hOmega^2$.
 Finally,
the data determining $g(\bx, t)$ is   the position
and velocity of the boundaries and so $g$ is periodic in $t$, 
and we get $\|\grad g(\bx, t)\|^2$ bounded above for all
time by its supremum
over one cycle, completing the time-independent
bound on the energy of $\bu(\bx, t)$.

\bibliography{iutampap}
\bigskip

\begin{tabular}{  p{6cm} }
  Philip Boyland \\
  Dept. of Mathematics\\
     University of Florida\\
     Little Hall\\
    Gainesville, FL 32605-8105 \\    
     \verb!boyland@math.ufl.edu! \\
\end{tabular}

\end{document}